\documentclass[11pt,a4paper]{article}

\usepackage{url}

\usepackage{latexsym,amssymb, amsthm}
\usepackage[centertags]{amsmath}
\usepackage{epsfig,pifont}
\usepackage[inactive]{srcltx}  

\usepackage{graphics,helvet,times,mathptm,epic,eepic}

\usepackage{color}

\usepackage{newlfont}

\usepackage{amsfonts,bbm}

\def\beq{\begin{equation}}
\def\eeq{\end{equation}}

\theoremstyle{remark}

\def\G1{\hbox{$\displaystyle{\mbox{\ding{172}}}$}}

\def\H2{\hbox{$\displaystyle{\mbox{\ding{173}}}$}}

\begin{document}

%

%
\title{Infinite Computations and the Generic Finite}

\newcommand{\nms}{\normalsize}
\author{  {   \bf Louis H. Kauffman\footnote{Louis H. Kauffman, Ph. D.,Professor at the University
of Illinois at Chicago, Chicago, Illinois, USA. }
    }\\ \\ [-2pt]
      \nms Department of Mathematics,\\[-4pt]
       \nms   University of Illinois at Chicago,\\[-4pt]
       \nms 851 South Morgan Street\\[-4pt]
       \nms Chicago, Illinois, 60607-7045\\ \\[-4pt]
       \nms http://www.math.uic.edu/$\sim$kauffman\\[-4pt]
         \nms {\tt  kauffman@uic.edu}
}

\date{}

\maketitle

%

\begin{abstract}
This paper introduces the concept of a {\it generic finite set} and points out that a consistent and 
significant interpretation of the {\it grossone}, $\G1$ notation of Yarolslav D. Sergeyev is that $\G1$ takes the role of a generic natural number. This means that $\G1$ is not itself a natural number, yet it can be treated as one and used in the generic expression of finite sets and finite formulas, giving a new power to algebra and algorithms that embody this usage. In this view, 
$$\mathbb{N} = \{1,2,3, \cdots , \G1-2,\G1-1,\G1\}$$ is not an infinite set, it is a symbolic structure representing a generic finite set.  
We further consider the concept of infinity in categories. An object $A$ in a given category
${\cal C}$ is {\it infinite relative to that category} if and only if there is a injection $J:A \longrightarrow A$ in ${\cal C}$ that is not a surjection. In the category of sets this recovers the usual notion of infinity. In other categories, an object may be non-infinite ({\it finite}) while its underlying set (if it has one) is infinite.
The computational  methodology due to Yarolslav D. Sergeyev  for executing
numerical calculations with infinities and infinitesimals is considered from this categorical point of view.
\end{abstract}



\section{Introduction}
We consider the grossone, $\G1,$ formalism of Yaroslav Sergeyev
\cite{Sergeyev,chaos,informatica,Menger,Calculus,Korea,patent,first,diff,Garro,Garro1,www,Rieman,Ludin}.
See also \cite{M,Z} for other discussions of this structure and its applications.
This paper begins in the stance that there are no completed infinite sets. 
We shall show that the grossone is then naturally interpreted
as a generic finite natural number and that the Sergeyev natural number construction
$$\mathbb{N} = \{1,2,3, \cdots , \G1-2,\G1-1,\G1\}$$ can be seen as a generic finite set.
By a {\it generic} finite set, I mean that $\mathbb{N}$ represents the properties of an arbitrarily chosen
finite set of integers taken from $1$ to a maximal integer represented by $\G1.$ This means that for 
any given finite realization of $\mathbb{\G1}$, $\G1$ will be larger than all the other integers in that realization. In this sense, we can say that $\G1$ is larger than any particular integer.
The grossone $\G1$ is a symbol representing the highest element in a generic finite set of natural numbers. Infinity is not the issue here. The issue is clarity of construction and the possibility of calculation with either limited or unlimited means.
\bigbreak

The original intent for the grossone and the set $$\mathbb{N} = \{1,2,3, \cdots , \G1-2,\G1-1,\G1\}$$ was that $\mathbb{N}$ should represent the infinity of the natural numbers in a way that makes the counting of infinity closer to computational issues than does the traditional Cantorian approach.  In Sergeyev's original approach, $\mathbb{N}$ is understood to be an infinite set, but one does not use the usual property of a Cantorian infinite set that an infinite set is in $1-1$ correspondence with a proper subset of itself. Rather, one regards special subsets such as the even integers as having their own measurement as parts of  $\mathbb{N}.$ For example, in the original approach it is fixed that grossone
$\G1$ is even, since it is postulated that the sets of even and odd numbers have the same number of elements, namely $\G1/2.$ Analogously, $\G1/3, \G1/4, \cdots , \G1/n$ are integral for finite $n.$
In our approach, which in this sense is a relaxation of the Sergeyev original approach, we interpret the grossone as a {\it generic} natural number and so it can be interpreted to be either even or odd.
We shall see, in Section $2$ how this has specific computational consequences in the summation of series. Both interpretations are useful. The generic integer intepretation that we use here has the advantage that the assumptions made are just those that can be applied to finite sets of natural numbers. This will become apparent as we continue with more theory and examples.
\bigbreak

The present paper is self-contained, and does not take on any axioms from the work of Sergeyev.
We use the notation $\G1$ as in the Sergeyev grossone, so that our work can be compared with his.
As we shall see, the result is that a different conceptual approach, via generic finite sets, has a very close relationship with the Sergeyev theory of an arithmetic of the infinite! Our theory can be seen as a relaxation of the Sergeyev theory where we no longer assume that $\G1$ is divisible by arbitrary natural numbers. As explained in Section 3, we adopt the 

\noindent {\it Transfer Principle: Any statement $P(\G1),$ using $\G1,$ is true if there is a natural number $N$ such that $P(n)$ is a true statement about finite natural numbers $n$ for all $n > N.$} 

\noindent This is the criterion determining the truth of statements about $\G1,$ creating a theory 
distinct from the Sergeyev theory, but sharing many of its notations and conceptual moves. See Section 3 for further discussion and, in particular, with a comparison with the Axioms of \cite{Z}.
\bigbreak

We reflect here many of these issues in the interpretation of $\mathbb{N}$ as a generic finite set. Our $\mathbb{N}$ is not a set but rather a symbolic structure that stands for any finite set.
As a symbolic structure $\mathbb{N}$ is inductively defined so that any of the following list of symbols denotes $\mathbb{N}:$
$$\mathbb{N} = \{1, \cdots , \G1\}$$
$$\mathbb{N} = \{1,2, \cdots , \G1-1,\G1\}$$
$$\mathbb{N} = \{1,2,3, \cdots , \G1-2,\G1-1,\G1\}$$
$$\mathbb{N} = \{1,2,3,4 \cdots ,\G1-3, \G1-2,\G1-1,\G1\}$$
$$\cdots$$
$$\mathbb{N} = \{1,2,3,\cdots ,n, \cdots ,\G1 - n+1,\cdots \G1-2,\G1-1,\G1\}$$
for any finite natural number $n.$
This means that there are infinitely many possible symbolic structures that indicate $\mathbb{N}.$
Each structure is finite as a symbolic structure. Depending upon the size of the set to which we shall refer, any and all of these symbols can be used for the reference. Any given positive integer  $n$ can occur, and does occur, in all such representatives past a certain point. Consequently,  we say that $n$ {\it belongs} to $\mathbb{N}$
for all natural numbers $n.$ And we say that $\G1$ {\it can be greater than} any given natural number $n.$ In this form we create a language in which speaking about $\mathbb{N}$ is very similar to speaking about an infinite set, but $\mathbb{N}$, being a symbol for a generic finite set, is not be interpreted as a Cantorian infinite set. The issues about subsets and $1-1$ correspondences do not arise. 
\bigbreak

The paper is organized into sections $2, 3$  and $4$ devoted to this point of view about genericity.
The style of writing for sections $2,3$ and $4$ is sometimes polemic, as it is written in the form of a speaker who is convinced that only finite sets should be allowed in mathematics. Infinite sets are seen to be taken care of by the concept of generic finite set. We show how to apply the grossone formalism by first thinking about how a computation could be written finitely, and then how it may appear in the limit where the grossone $\G1$ is itself regarded as an infiinite number. Except for certain situations where we can consider limits as $\G1$ becomes infinite, we do not use the original intepretations of the grossone as an infinite quantity. An example of a situation where one can shift from one interpretation to the other is in a summation of the form $S = \Sigma_{n=1}^{\G1} F(n)$ where $F$ is a function defined on the natural numbers. For example, consider $$\Sigma_{n=1}^{\G1} 2^{n} = 2^{\G1 + 1} - 1.$$
We can regard this as standing for all specific formulas of the type 
$$\Sigma_{n=1}^{N} 2^{n} = 2^{N + 1} - 1$$ where $N$ is any natural number, {\it or} we can regard it as an generic  summation with an generic result of $ 2^{\G1 + 1} - 1.$ The generic result can be regarded in as formally similar to the infinite result that comes from taking $\G1$ infinite. Both ways of thinking about the answer tell us how the summation behaves when the number of summands is $\G1$ and $\G1$ very large. Other examples will be discussed in the body of the paper. 

\bigbreak

Secion $5$ is written from another point of view. In this section we take a categorical and relative point of view about ``being infinite." We accept that in bare set theory a set is infinite if it is in $1-1$ correspondence with a proper subset of itself. But in another category we ask that this $1-1$ correspondence be an injection in that category that is not a surjection. This means, for example, that a circle or a sphere in the topological category is not infinite (hence finite by definition) since there is no homeomorphism of a sphere to a proper subset of itself. We examine the model $\mathbb{N}$ for the extended natural numbers from this point of view, and show that it can be construed as set theoretically infinite, but topologically finite. We hope that these two points of view for interpretation will enrich the intriguing subject of the grossone and it uses.
\bigbreak

\section{The grossone, ``Infinite" Series and The Three Dots}
In Yaroslav Sergeyev's theory of  numeration, he considers a completion of the natural numbers
$${\cal N} = \{1,2,3,\cdots \}$$ to a set that contains an infinite number $\G1,$ referred to as
{\it grossone}. The completion is denoted by the notation 
$\mathbb{N}$ in the form
 \beq
\mathbb{N} = \{ 1,2,  \hspace{3mm} \ldots  \hspace{3mm}  \ldots \hspace{3mm} \G1-2,
\hspace{2mm}\G1-1, \hspace{2mm} \G1 \}. 
       \eeq
The grossone, $\G1$, behaves ``just like" a very large integer, so that $\mathbb{N}$ as a set is 
conceived to have all of its properties analagous to those of a finite set of integers such as
$$Set(N) = \{ 1,2,  \hspace{3mm} \ldots  \hspace{3mm}  \ldots \hspace{3mm} N-2,
\hspace{2mm} N-1, \hspace{2mm} N \}.$$
While Segeyev does not quite explicitly say that $\mathbb{N}$ and $Set(N)$ are logical twins just so long as $N$ is very very large, this is the basic idea that we explore in this paper. For example, it should not be the case that 
$\mathbb{N}$  should appear infinite in the sense of having a $1-1$ correspondence with a proper subset of itself. Some of the usual attempts obviously fail. For example if we try to map  $\mathbb{N}$  to itself by the map $f(k) = k +1$ then $f(\G1) = \G1 + 1$ and so we need a larger domain.
This mirrors the problem for the finite set into the infinite set $\mathbb{N}.$
\bigbreak

But what set is this $\mathbb{N}$? How should we interpret its two appearances of the ``three dots"?
In the case of ${\cal N}$ the three dots refer to the Peano axioms for the (usual) natural numbers, assuring us that given a natural number $n$ then there is a successor to that number indicated by
$n + 1,$ and that $n+1$ is never equal to $n.$ The principle of mathematical induction is then used to characterize the set ${\cal N}.$ If $M$ is any subset of  ${\cal N}$ containing $1$ and having the property that $n \in M$ implies that $(n+1) \in M,$ then $M = {\cal N}.$
\bigbreak

At first glance one is inclined to guess that $\mathbb{N}$   consists in two copies of ${\cal N},$
one ascending $$\{1,2,3, \cdots \},$$ and one descending $$\{ \cdots \G1-3,\G1-2, \G1-1, \G1 \}.$$
If we then take the union 
$$\{1,2,3, \cdots  \cdots \G1-3,\G1-2, \G1-1, \G1 \},$$ is this  $\mathbb{N}?$ I submit that this {\it is not}
the desired $\mathbb{N}.$ For one thing, if one starts counting upward $1,2, \cdots$ one never leaves the left half of this union. And if one starts counting down, one never leaves the right half of the union.
This is not in analogy to the large finite set, where either counting down or counting up will cover all of the territory. Also, there appears to be an injection of the union to itself that is not a surjection. Just add one to every element of the left half and subtract one from every element of the right half. All in all, we must search for a different model for $\mathbb{N},$ or a different category for it to call its home.
\bigbreak

I believe that the simplest interpretation for $\mathbb{N}$ is that it is a {\it generic finite set.}
This means that we interpret $\G1$ as a generic ``large" natural number. $\G1$ is not any specific natural number unless we want to take it as such. Then since $\mathbb{N}$ is a generically finite set, it has no $1-1$ correspondence with a proper subset of itself. In fact, $\mathbb{N}$  is not a set at all.
It is a symbolic construct that represents the {\it form} of a finite set. It is not a set just as an algebraic variable $x$, standing for a number, is not a number. As such,
 $\mathbb{N}$ has no members, nor does it have any subsets. As a symbolic construct however, 
 $\mathbb{N}$ has some nice features. We all agree on the equalities
 $$\mathbb{N} = \{1, \cdots ,\G1-1,\G1\},$$
 $$ = \{1,2, \cdots , \G1-2,\G1-1,\G1\},$$
 $$ = \{1,2,3, \cdots , \G1-3, \G1-2,\G1-1,\G1\},$$
 $$= \cdots$$
 $$= \{1,2,3, \cdots, n, \cdots, \G1-n, \G1-(n-1), \cdots,  \G1-3, \G1-2,\G1-1,\G1\},$$
 for any {\it specific} natural number $n.$ This is the nature of  $\mathbb{N}$ as a symbolic construct.
 In this sense it gives the appearance of acting like an infinite set. This property of the notations that we use is an inheritance from the set theory notation where exactly the same symbolic phenomenon seems to indicate infinity where one only has the form of a series of numbers and the Peano property of knowing that for a given $n$ there is an $n+1.$ In the usual notation for the set of natural numbers, we have $${\cal N} = \{1,2,3, \cdots, n, \cdots\}$$ for any specific natural number $n.$ Thus we can also regard ${\cal N}$ as a generic symbol, but it is a generic symbol that is incomplete if it is to be viewed as 
 representing a finite set. The symbolic construction for $\mathbb{N}$  allows us to start anew and eliminate the notion of a completed infinity from our work.
 \bigbreak

If we interpret $\mathbb{N}$  as some particular finite set, then it has $N$ members where $N$ is some specific finite natural number. It can have either an odd or an even number of members. It is subject to permutations just as any finite set is so subject. 
The availability of such interpretations makes our viewpoint different from the original interpretation of the grossone. In the original interpretation, $\G1$ is always infinite, and it is divisible by any natural number. In our interpretation $\G1$ as a symbolic entity is not a finite natural number. It stands for any top value in some finite set, but not any particular value. For this reason we can do arithmetic with $\G1$ and say that $\G1$ is taken to be greater than $n$ for any standard natural number $n$ and that $\G1 -1 < \G1 < \G1 + 1 < 2^{\G1}.$ 
\bigbreak

{\it In the generic context, for a  statement about $\G1$ to be true, it must be true for all sufficiently large substitutions of natural numbers for $\G1.$} For example $\G1 > 2^{100}$ is true because $N > 2^{100}$ is true for all natural numbers $N$ that are larger than $2^{100}.$ Some statements are true for all numbers. Thus $2^{\G1} > \G1$ is true since $2^{N} > N$ for all natural numbers $N.$ In this way, we can use the grossone as an infinite number when we wish.
\bigbreak

\noindent Series of the form $$\Sigma_{k=1}^{\G1} a_{k},$$ where $a_{k}$ are any real or complex numbers or algebraic expressions, are well defined. For example, we can write
$$S = S(\G1) =  1 + x + x^2 + \cdots + x^{\G1},$$ and this is a generic finite sum in the form of a geometric series.
We can operate algebraically on such sums. For example,
$$xS = x(1 + x + x^2 + \cdots + x^{\G1}) = x + x^2 + \cdots + x^{\G1} + x^{\G1 + 1},$$
whence
$$(1-x)S = 1 - x^{\G1 + 1},$$
and so, for $x \ne 1$
$$S(\G1) = \frac{1 - x^{\G1 + 1}}{1-x}.$$
\bigbreak

Note that we can replace $\G1$ in this formula. For example,
$S(2 \G1) = \frac{1 - x^{2 \G1 + 1}}{1-x}.$ It is our intent to keep the interpretation that $\G1$ is 
a generic natural number, so in this sense $S(\G1)$ can be regarded as specializing to the function
$S(n) =  \frac{1 - x^{n + 1}}{1-x},$ for natural numbers $n.$ But $S(\G1)$ is also the generic expression
for this function, and as such we can think of it as expressing cases where $n$ is taken to infinity.
\bigbreak

At this point we can note the effects of taking $\G1$ even or odd, or the effect of taking a limit as 
$\G1$ becomes arbitrarily large. For example if the absolute value of $x$ is less than 1, then
$x^{\G1 + 1}$ becomes arbitrarily small as $\G1$ becomes arbitrarily large. Generically, we can
regard $\G1$ as an infinite number if we wish. In that case we have that $S$ is infinitesimally close to
$ \frac{1}{1-x},$ but the expression $S = \frac{1 - x^{\G1 + 1}}{1-x}$ is a more accurate rendition of the actual situation. Working with generic formalism allows us to dispense with limits in many cases and adds detail that can be sometimes ignored in the usual category of working with limits.
\bigbreak

Note in the last example, that if we take $x=-1,$ then 
$$S =  \frac{1 - (-1)^{\G1 + 1}}{2}=  \frac{1 + (-1)^{\G1}}{2}.$$
If we take $\G1$ even, then $S = 1.$ If we take $\G1$ to be odd, thn $S=0.$
We see that $(-1)^{\G1}$ is a well-defined symbolic value that can be either positive or negative, depending on its instantiation as a number. This is in direct contrast to Sergeyev's usage for the
grossone, where $\G1$ is taken to be even and so $S$ takes only the value $1.$ Our interpretation reflects the fact that the corresponding finite series oscillate between $0$ and $1.$
\bigbreak

Here is another example. Let
$$E(x) = (1 + \frac{x}{\G1})^{\G1}.$$
This is of course a rendition of $(1 +x/N)^{N}$ for large generic $N.$ We can apply the binomial theorem
to conclude that 
$$E(x) = 1 + \frac{\G1}{1!}(\frac{x}{\G1}) +\frac{\G1(\G1 - 1)}{2!}(\frac{x}{\G1})^{2} + \cdots +   \frac{\G1(\G1 - 1)(\G1 -2) \cdots (1)}{\G1 !}(\frac{x}{\G1})^{\G1}.$$ Thus
$$E(x) = 1 + \frac{x}{1!}  + \frac{\G1(\G1 - 1)}{\G1^{2}} \frac{x^2}{2!} + \frac{\G1(\G1 - 1)(\G1-2)}{\G1^{3}}\frac{x^3}{3!} +  \cdots +  (\frac{x}{\G1})^{\G1}$$
From this expression we can see the limit structure that leads to the usual series formula for 
$$e^{x} =  1 + \frac{x}{1!}  + \frac{x^2}{2!} + \frac{x^3}{3!} +  \cdots ,$$
and we see the luxury in the exact formula for $E(x).$ By not taking the limit and examining the exact 
structure of the formula for large generic $\G1,$ we see more and can write down more exact approximations for specific values $N$ substituted for the generic $\G1.$ 
\bigbreak

Another example will show even more. Let $$\mathbb{P} = \{ p_{1}, p_{2}, \cdots , p_{\G1-1}, p_{\G1} \}$$
denote a list of all the prime numbers up to a generic natural number $\G1.$ Note that $p_{\G1}$ 
denotes a {\it generic prime number}. This can be compared with Sergeyev's concept of 
infinite prime numbers \cite{Rieman}.
\bigbreak

Define an analogue of the 
Riemann Zeta function via
$$Z(s) = [ \frac{1}{1- \frac{1}{p_{1}^{s}}}] [\frac{1}{1- \frac{1}{p_{2}^{s}}}] \cdots 
[\frac{1}{1- \frac{1}{p_{\G1}^{s}}}] = \Sigma_{n \in N(\G1)} \frac{1}{n^{s}}.$$ 
Here $N(\G1)$ denotes all
natural numbers that can be constructed as products of prime powers from the set of primes 
$\mathbb{P}.$ Note that there is no limit to the size of elements of $N(\G1).$ If we wish to keep bounds on that, we will have to introduce further notation. Even though $Z(s)$ is a finite product, it produces natural numbers of arbitrarily high size. Thus a better finite zeta function will be given by 
$$ZZ(s) = \Pi_{k=1}^{\G1}(1 + p_{k}^{-s} + p_{k}^{-2s} + \cdots + p_{k}^{- \G1 s}).$$
Then 
$$ZZ(s) = \Sigma_{n \in N(\G1,\G1)} \frac{1}{n^{s}}$$ where $N(\G1,\G1)$ denotes the generically finite set of natural numbers whose prime factorization is from the first $\G1$ primes, and whose prime powers are no higher than $\G1.$ For a particular finite instantiation of $\G1$ this zeta function, being a finite sum, can be computed for any complex number $s.$ This truncation of the usual limit version of the zeta function
allows new computational investigations of its properties. This example can be compared with the work
in \cite{Rieman} and \cite{Z}. One possible advantage of our approach is that it is essentially a finite approach. We are suggesting that it is useful to examine the Rieman zeta function, defined just for a finite collection of prime numbers, and then to examine how these finitized zeta functions behave as the number of primes used gets larger and larger. We are not concerned here with extra fine structure of the zeta function at its limit, but rather with extra structure visible in the finite approximations to it. This material will be explored further in a separate paper.
\bigbreak

But the reader may ask, does not the construction of $\mathbb{N}$
connote an infinite set if we regard the grossone $\G1$ as standing for an infinite evaluation, larger than any natural number? Well, dear reader, what will you have? Should we not be able, in this putative infinity, to count down from $\G1$ by successive subtractions to some number that we found by counting up from $1$ by successive additions? This is true of any finite set. Yet if this is true then it must be that for some natural numbers $n,m$ we have $$n = \G1 -m$$ and so we find that $$n+m = \G1$$ and so the set
$\mathbb{N}$ would be neither generic, nor infinite, but simply one of the multitude of finite sets, in fact, the one of the form$$\{1,2,3, \cdots , n+m-2, n+m-1, n+m\}.$$ Once again, 
$$\mathbb{N} = \{1,2,3, \cdots, \G1-2, \G1-1,\G1\}$$ is not a set at all, but the form of a generic finite set.
It is a symbol for the finite set structure, not a set at all. But this gives us freedom to regard this symbolic structure as something like a set, and in it, we indeed cannot countdown from the grossone to a finite natural number.
\bigbreak

\section{Infinity and the Generic Finite - The Transfer Principle}
In many aspects of mathematics there is no need for any infinite set. It is sufficient to have the concept of a generic finite set. For this we need the notion of a generic natural number. Thus one speaks of the 
set of natural numbers from $1$ to $n$, $\{1,2,3, \cdots n-2,n-1,n \}$, and one usually conceives this as 
referring to some, as yet unspecified number $n.$ One writes a formula such as 
$$1 + 2 + \cdots + n = n(n+1)/2,$$ and it is regarded as true for any specific number $n.$ The concept of a generic $n$ requires a shift of attention, but no actual change in the formalism of handling finite sets and finite series. To write
$$1 + 2 + \cdots + \G1 = \G1(\G1+1)/2$$ is conceptually different from the above formula with $n$.
In the formula with $n$ we are referring to {\it some specific natural number $n$}. In the formula with grossone, we are indicating a {\it form} that is true when $\G1$ is replaced by a specific integer $n,$ and we are also indicating the behaviour of a corresponding limit or infinite sum.
\bigbreak

\noindent {\bf The Transfer Principle.}
Once we take a notation for a generic natural number as with Sergeyev's grossone,
we are lead to use the concept of genericity. In this view $\G1$ is symbolic, not infinite, but can be regarded as indefinitely large. We can regard it as larger than any given number that is named. This is a way of thinking instantiated by our {\it Transfer Principle: Any statement $P(\G1),$ using $\G1,$ is true if there is a natural number $N$ such that $P(n)$ is a true statement about finite natural numbers $n$ for all $n > N.$} This principle provides a transfer statement that allows us to apply the grossone in many particular situations. The principle and its consequences are distinct from the Sergeyev use of the grossone, and constitute a relaxation of that usage.
\bigbreak

 For example, if we are working with a computer that is limited to a specific size of natural number, then from the outside, as theorists, we can easily say that $\G1$ will be greater than any number allowed in that computational domain. We can think of $\mathbb{N}= \{1,2, 3,\cdots , \G1-2, \G1-1, \G1\}$ as a set that is larger than any specific finite set we care to name, but it is still generically finite and does not partake of the Cantorian property of being in 
$1-1$ correspondence with a proper subset of itself. It is like a finite set, but it is not any particular finite set.
\bigbreak

The concept of a generic set is different from the concept of a set in the same way that the variable $x$ is different from a specific number in elementary algebra. A generic set such as $\mathbb{N}$ does not have a cardinality in the sense of Cantor. It is not in $1-1$ correspondence with any specific finite set, but if a specific natural number $n$ is given for $\G1$, then the resulting set is of cardinality $n.$
This is precisely analogous to the situation with an algebraic $x$ that in itself has no numerical value, but any subsititution of a number for $x$ results in the specific value of that number. Just as we endow algebraic expressions with the same properties as the numbers that they abstract, we endow generic sets with the properties of the finite sets which they stand for. 
\bigbreak

\noindent {\bf The Zhigljavsky Axioms.}
The Transfer Principle can be compared with axiomatizations for the Sergeyev system. For example,
in \cite{Z}, Zhigljavsky takes the following set of axioms:
\begin{enumerate}
\item Grossone, $\G1$ is the largest natural number so that 
$\mathbb{N} = \{1,2,3, \cdots ,\G1-1, \G1\}.$
\item The grossone $\G1$ is divisible by any finite natural number $n.$
\item If a certain statement is valid for all $n$ large enough, then this statement is also valid for
$\G1.$
\item Assume we have a numerical sequence $\{b_{n}\}_{n=1}^{\infty}$ which converges to $0$ as
$n \longrightarrow \infty.$ Then this sequence has a grossone-based representative 
$\{ b_{1},b_{2}, \cdots b_{\G1}\}$ whose last term $b_{\G1}$ is necessarily an infinitesimal quantity.
\end{enumerate}
We see that Axiom 1 corresponds to our formalization of $\mathbb{N}$ as a generic set of consective natural numbers. But we do not use formally the concept of the ``set of all natural numbers". We do not assume Axiom 2. Our transfer principle is essentially the {\it same} as Axiom 3. We agree with 
Axiom 4 in the sense that one has the generic set $\{ b_{1},b_{2}, \cdots b_{\G1}\},$ and the symbol
$b_{\G1}$ can be used to stand for a number whose absolute value is less than that of any given finite non-zero real number. (In fact it follows from the Transfer Principle that $|b_{\G1}|< 1/n$ for any positive
natural number $n.$)  This means that in both our system and in the  Zhigljavsky system, infinitesimals naturally arise. How they are treated may depend upon the further adoption of a theory of extended real analysis and is left open for specific applications.
\bigbreak

Our point in this paper has been that the 
symbolic constructions  of Yaroslav Sergeyev can be regarded as generic finite sets. Does this preclude an infinite interpretation? In the case of a series such as 
$$S = 1 + x + x^2 + \cdots + x^{\G1}=  \frac{1 - x^{\G1 + 1}}{1-x}$$ ($x$ not equal to  $1$), we would like to give the interpretation that if $\G1$ is an  infinite integer, and $|x| < 1$, then S is infinitesimally close to $1/(1-x).$ The problem with this is the same as the corresponding problem in the calculus. An integral is a limit of finite summations. We are led to imagine (by the Leibniz \cite{Leibniz,Newton} notation for example) that the integral $$\int_{-\infty}^{\infty} e^{-x^2 /2}dx$$ is an uncountably infinite sum of infinitesimal terms of the 
form $e^{-x^2 /2}dx.$ It takes the language of non-standard analysis \cite{Robinson} to make formal sense out of this statement. We can, using the grossone formalism, go in the other direction and articulate that integral as a generic finite sum. In order to do this we have to make a choice of method of integration and then write the formula in generic fashion. For example, consider  $$\int_{-\infty}^{\infty} f(x)dx.$$ We use Riemann integration and choose an interval from $-n$ to $n$ on the real line. We divide the interval $[-n,n]$ up into sub-intervals of length $1/n.$ This gives the partition
$$\{ -n^2/n, (-n^2 + 1)/n, (-n^2 + 2)/n, \cdots, -1/n, 0 , 1/n, 2/n, \cdots, (n^2 - 1)/n, n^2/n \}$$ of the interval
$[-n, n].$ We then have that, when the integral converges, 
$$\int_{-\infty}^{\infty} f(x)dx = limit_{n \longrightarrow \infty}\, \Sigma_{k= -n^2}^{k= n^2} \frac{f(k/n)}{n}.$$
Thus the following finite generic expression stands for this integral.
$$\Sigma_{k = -\G1^2}^{\G1^2} f(\frac{k}{\G1})\frac{1}{\G1}.$$ We can even write
$$\int_{-\infty}^{\infty} e^{-x^2 /2}dx = lim \,\, \Sigma_{k = -\G1^2}^{\G1^2} E(\frac{k^2}{2 \G1})\frac{1}{\G1}$$ where $E(x)$ is our finite version of the exponential function from the previous section. Here 
$lim \,\, F(\G1) = lim_{n \longrightarrow \infty}\,\, F(n)$ where the limit is the classical limit and $n$ runs over the classical natural numbers. We cannot write a limit as $\G1$ approaches something since $\G1$ is generic and does not approach anything other than itself.
\bigbreak

Let us write $a \doteq b$ to mean that $a$ and $b$ differ by an infinitesimal amount (in the sense of the discussion above) or, equivalently, that $lim \,  a = lim \,  b.$ Thus 
$$ \frac{1 - x^{\G1 + 1}}{1-x} \doteq  \frac{1}{1-x}$$ when $|x| < 1.$
Note that $a \doteq b$ has the following properties, making it an equivalence relation (the first three properties) and more.
\begin{enumerate}
\item $a  \doteq  a$ for any $a.$
\item If $a  \doteq b$, then $b  \doteq a.$
\item If $a  \doteq b$ and $b  \doteq c$ then $a  \doteq c.$
\item if $a \doteq b$ and $c \doteq d$ then $a + b \doteq c + d$ and $ac \doteq bd.$
\item If $a  \doteq b$ then $a^{1/n}  \doteq b^{1/n}$ where $x^{1/n}$ denotes the prinicipal $n$-th root of the number $x$ (working over the complex numbers). This implies that $$a^{1/\G1}  \doteq b^{1/\G1}$$
in appropriate cases where one checks the limiting behaviour.
\item If $1/\G1 \doteq \epsilon$ where $\epsilon$ is infinitesimal, then it may happen that 
$1 \doteq \G1 \epsilon.$ In this case one must examine the size of $\epsilon$ in relation to the size of 
$1/\G1.$
\end{enumerate}
\bigbreak

\noindent Here is an example. In the derivation that follows we will use all of the above properties of 
$a \doteq b,$ and we leave it for the reader to check that each of the steps is valid.
We have the classical formula $$e^{i \theta} = cos(\theta) + i sin(\theta).$$
This is often seen as a consequence of the series formula for $e^{x}.$ We can write the equation
$$(1 + \frac{i\, \theta}{\G1})^{\G1} \doteq cos(\theta) + i \, sin(\theta).$$ In particular, the well-known formula $e^{i \pi} = -1$ becomes the following generic limit formula
$$(1 + \frac{i\, \pi}{\G1})^{\G1} \doteq -1.$$ 
In fact, we can, since this formula refers to properties of all formulas where $\G1$ is replaced by a specific integer, take roots and solve for $\pi$ as follows.
$$(1 + \frac{i\, \pi}{\G1})^{\G1} \doteq -1$$ 
$$1 + \frac{i\, \pi}{\G1} \doteq (-1)^{\frac{1}{\G1}}$$ 
$$\frac{i\, \pi}{\G1} \doteq (-1)^{\frac{1}{\G1}} -1$$ 
$$i\, \pi \doteq \G1 ((-1)^{\frac{1}{\G1}} -1)$$ 
$$\pi \doteq \frac{\G1 ((-1)^{\frac{1}{\G1}} -1)}{i}$$ 
One can verify that this last limit formula is indeed a correct limit formula for $\pi.$
Well-known limit formulas appear from this formula when we replace $\G1$ by $2^{\G1}$ in it. We then have 
$$\pi \doteq \frac{2^{\G1} ((-1)^{\frac{1}{2^{\G1}}} -1)}{i}$$ 
In this version the finite versions are of the form
$$ \frac{2^{N} ((-1)^{\frac{1}{2^{N}}} -1)}{i}$$ 
and the limit formula is
$$\pi = lim_{N \longrightarrow \infty} \frac{2^{N} ((-1)^{\frac{1}{2^{N}}} -1)}{i}$$ 
$$\pi = lim_{N \longrightarrow \infty} 2^{N} Imag((-1)^{\frac{1}{2^{N}}}),$$  where $Imag$ denotes the imaginary part of the given complex number.
Here we have successive square roots of $-1$ and we can use the formula
$$\sqrt{a + bi} = \sqrt{(1+a)/2} + i \sqrt{(1-a)/2}$$ when $a^2 + b^2 = 1.$ 
From this it is easy to derive the famous formula of $Vi\grave{e}te:$
$$\pi = \frac{2}{ \sqrt{\frac{1}{2}} \sqrt{\frac{1}{2} +\frac{1}{2} \sqrt{\frac{1}{2}} }   \sqrt{\frac{1}{2} + \frac{1}{2}  \sqrt{\frac{1}{2} +\frac{1}{2} \sqrt{\frac{1}{2}} }  }        \cdots    }$$
The grossone notation and the notion of generic finite sets allows us to write this derivation in a concise and precise manner.
\bigbreak

\section{Limits, Ordinals and Further Relaxations}
In articulating the notion of the generic finite we have examined an interpretation of the grossone extension of the natural numbers of Yaroslav Sergeyev as a generic finite set. This works because the grossone extension puts a cap, the grossone $\G1$, at the top of its model of the natural numbers and so formally resembles our way of thinking about a finite set such as $\{ 1,2,3 \}$ where there is a least element, a greatest element, and a size for the set in the numerical sense of counting or ordinality.
By working with this correspondence we may end up with evaluations of size that are dfferent  than the Sergeyev theory.
As case in point is the finite sets of even natural numbers.
$$\{2 \},$$
$$\{2,4 \},$$
$$\{2,4,6 \},$$
$$\cdots$$
$$\{2,4,6,8,\cdots, 2n\}.$$
The {\it generic set} of even natural numbers is 
$$E = \{2,4,6,8, \cdots 2\G1 \}.$$
This notation is different from the way Sergeyev would denote the even natural numbers, since he would place them ``inside" the set $$\mathbb{N} = \{1,2,3, \cdots , \G1-2,\G1-1,\G1\}.$$ As such there could be no appearance of $2 \G1$ in the Sergeyev representation of the even numbers.
In Sergeyev's language the {\it size} of the set of even natural numbers is  $\G1/2,$ since it comprises half of the natural numbers. We can discuss {\it size}  in our way by the transfer principle of the previous section. The size of a specific finite set of even numbers 
$$\{2,4,6, \cdots 2n \}$$ is the number of elements of the set, which is $n$ in this case. If the size for the specializations (replacing $\G1$ by specific natural numbers $n$) of a generic set $X$ has the form $F(n)$ for a function $F$ of natural numbers, then we say that {\it $X$ has size $F(\G1).$} Thus we assign 
$\G1$ as the size for $E$ above. This moves our theory in the direction of Cantorian counting for sets even though we have not invoked a notion of $1-1$ correspondence for infinite sets.
The point of this example is to show that there are are natural differences between the generic finite set approach and the Sergeyev approach to the grossone.
\bigbreak

There are other relationships of the generic finite set concept and standard set theory.
For example, in standard set theory we take as representative ordinals
$$0 = \{ \, \},$$
$$1 = \{0\},$$
$$2 = \{0,1\},$$
$$3= \{0,1,2\},$$
and so the {\it generic finite ordinal} is 
$${\cal G} = \{0,1,\cdots , \G1\}.$$ This should be held in contrast to the first infinite ordinal
$$\omega = \{ 0,1,2,\cdots \}.$$
The generic finite ordinal has $\G1 + 1$ members and so we could name it 
$$\G1 + 1 = \{0,1,\cdots , \G1\}$$
in the tradition of making ordinals. However, ordinals made in this generic fashion are not well-ordered since there is no end to the descending
sequence $\G1, \G1-1, \G1-2 , .\cdots.$ This means that a theory of ordinals based on generic finite sets will have a character of its own. This will be the subject of a separate paper.
\bigbreak

\section{Grossone and Infinity Relative to a Category}
In this section, we take a different approach to the grossone. We assume here the existence of infinite sets and the usual terminology of point set topology. With that we discuss, by using categories, relative notions of infinity. After all, a circle is not homeomorphic to any proper subset of itself. Therefore the point set for the circle, uncountable in pure set theory, is {\it finite} in the category of topological spaces!
We formalize this notion below and indicate how it can be interfaced with the grossone.
\bigbreak

We recall that a {\it category} ${\cal C}$ is a collection of {\it objects} and {\it morphisms} where a morphism is  associated to two objects and is usually written as $f:A \longrightarrow B$ where 
$A$ and $B$ are objects. Note that the morphism $f$ provides a directed arrow from $A$ to $B.$
Without further axioms the concept of a category is the same as the concept of a directed multi-graph.
The axioms for a category are as follows:

\begin{enumerate}
\item Given morphisms $f:A \longrightarrow B$ and $g:B \longrightarrow C,$ there is a well-defined
morpism called the {\it composition of $f$ and $g$} and denoted 
$$g \circ f: A \longrightarrow B.$$ The object $A$ is called the {\it domain} of $f$ and the object $B$ is
called the {\it codomain} or {\it range} of $f.$

\item Every object $A$ has a unique {\it identity morphism} $$1_{A}: A \longrightarrow A$$ such that
for any $f:A \longrightarrow B$, $f \circ 1_{A} = f$ and for any $g:B \longrightarrow A$, $1_{A} \circ g = g.$

\item If $f:A \longrightarrow B,$  $g:B \longrightarrow C$ and  $h:C \longrightarrow D,$ then
$$h\circ (g \circ f) = (h\circ g) \circ f.$$ Thus composition of morphisms is associative.
\end{enumerate}

If $C$ and $C'$ are categories, then we say that a {\it functor} from $C$ to $C'$, denoted
$F:C \longrightarrow C'$, is a function that takes objects to objects and morphisms to morphisms
such that if $f:A \longrightarrow B$ is a morphisim in $C$, then $F(f):F(A) \longrightarrow F(B)$ is a morphism in $C'.$ Furthermore, we require of a functor that identity morphisms are carried to identity morphisms, and that compositions are taken to compositions in the sense that $F(f \circ g) = F(f)\circ F(g)$ for all compositions $f \circ g$ in $C.$
\bigbreak

In this paper we will use categories whose objects are sets and whose morphisms are maps of these sets with whatever extra structure is demanded by that category. Then there is a forgetful functor
$FG: {\cal C} \longrightarrow {\cal Sets}$ obtained by just taking the objects as sets and the morphisms as maps of sets, ignoring the extra structure. In such categories, a morphism $f:A \longrightarrow B$
is said to be {\it injective} if $FG(f)$ is injective in ${\cal Sets},$ and $f$ is said to be 
{\it surjective} if $FG(f)$ is surjective in ${\cal Sets}.$
\bigbreak

In ${\cal Sets}$ one says that a set $A$ is {\it infinite} if there exists an injection $i:A \longrightarrow A$
that is not a surjection. If every injection of $A$ to itself is a surjection, we say that the the set $A$ is
finite. We relativize this notion to other categories. If ${\cal C}$ is a set-based category, we say that 
an object $A$ of ${\cal C}$ is {\it finite} (in ${\cal C}$)  if every injection of $A$ to itself {\it in the category ${\cal C}$ } is
surjective. For example consider the category ${\cal Top}$ of topological spaces. We see at once that 
the circle  $$S^{1} = \{(x,y)|x^2 + y^2  = 1\}$$ where $x,y$ are real numbers, is finite in this category since the circle is not homeomorphic to any proper subset of itself. Thus, while the circle  consists in infinitely many points when looked at under the forgetful functor to set theory, in the topological category the circle is finite. We shall see that this point of view on finite and infinite is very useful in sorting out how we deal with mathematical objects in many situtations.
\bigbreak

\begin{figure}
     \begin{center}
     \includegraphics[width=8cm]{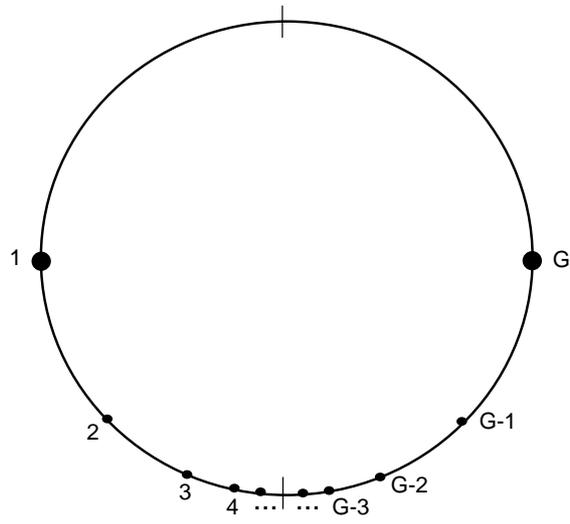}
     \caption{Grossone Circle}
     \label{gcircle}
\end{center}
\end{figure}

\begin{figure}
     \begin{center}
     \includegraphics[width=8cm]{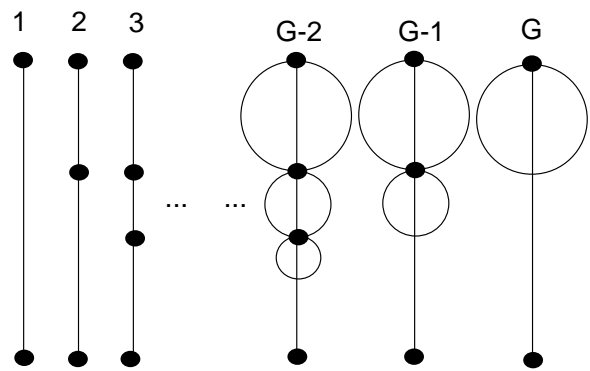}
     \caption{Grossone Graphs}
     \label{graphs}
\end{center}
\end{figure}

We now combine this topological point of view on finiteness with the grossone.
In Figure~\ref{gcircle} we show an embedding of a union of two sets of points in the form
$$\mathbb{N} = \{1,2,3,3,\cdots \,\, \cdots, \G1-3,\G1-2,\G1-1,\G1 \}.$$ (In the Figure, $\G1$ is denoted by $G.$)
We regard $\mathbb{N}$ as embedded
in a circle with two points removed. The two vertical marks on the circle in the figure denote the removed points. The circle with two points removed partakes of the subspace topology from the Euclidean plane in which it is embedded. We work in the category of orientation preserving homeomorphisms of the deleted circle ${\cal S}$
that map the the intervals $[i,i+1]$ and $[\G1-i,\G1-i+1]$ to themselves. In the categorical sense, ${\cal S}$ is finite and there can be no such homeomorphisms that take $\mathbb{N}$ to a proper subset of itself.
In fact, in this model, every such homeomorphism is the identity map when restricted to $\mathbb{N}.$
\bigbreak

In Figure~\ref{graphs} we give another example of how to topologically make an infinite set finite. 
We have labeled a set of graphs with the ``elements" of $\mathbb{N}.$ None of these graphs are homeomorphic (we take the nodes of the graphs to be disks and the edges to be topological intervals)and so a homeomorphism of the entire collection must take each graph to itself. There is no topological injection of the collection of graphs to a subcollection of itself.
\bigbreak

We produce this model to suggest some ways to view $\mathbb{N}$ as infinite and yet finite without paradox. I am sure that other models will emerge in relation to applications of these ideas.
\bigbreak

\bigskip


%


\normalsize

\end{document}